# A Fast Solver for Tridiagonal Toeplitz Systems with Multiple Right-Hand Sides

A Preprint

**Shahin Hasanbeigi**
Department of Applied Mathematics
Tarbiat Modares University
Tehran, Iran
shahin77hb@gmail.com

June 2024

**Abstract**

In this work, we introduce a novel approach for solving tridiagonal Toeplitz systems with multiple right-hand sides. Tridiagonal Toeplitz matrices are a prevalent structured matrix type that appear in numerous applications and fields, including image processing, partial differential equations (PDEs), fast Fourier transform (FFT), and others. In this work, we extend and utilize the fast solver introduced in [7] to solve systems of the form $(I \otimes T)x = b$. At the conclusion, numerical examples demonstrate the efficacy and precision of the introduced algorithm.

***Keywords*** Tridiagonal Toeplitz Systems · Multiple Right-Hand Sides · Fast Solver · Sherman-Morrison-Woodbury Formula

## 1 Introduction

Consider the following nonsingular linear system of equation

$$TX = B \tag{1}$$

where T is an $n \times n$ Tridiagonal Toeplitz matrix of the form

$$T = \begin{bmatrix} t_1 & t_3 & & & & \\ t_2 & t_1 & t_3 & & & \\ & \ddots & \ddots & \ddots & & \\ & & \ddots & \ddots & \ddots & \\ & & & t_2 & t_1 & t_3 \\ & & & & t_2 & t_1 \end{bmatrix} \tag{2}$$

$X$ is $n \times m$ unknown matrix, and $B$ is $n \times m$ right-hand side matrix.

Systems of equations of this nature arise in a multitude of applications and fields, including but not limited to image and signal processing, discretizing partial differential equations, and so forth. In general, there are two principal approaches to solving any linear system of equations: direct and iterative methods. Iterative methods are primarily comprised of Krylov subspace methods, such as the CG method for Hermitian positive definite matrices and the GMRES method for non-Hermitian positive definite matrices, as well as their variants, as cited in Chan et al.[1, 5, 3], these include classical splitting methods such as Gauss-Sidel, SOR, and others. In contrast, direct methods, as described in [2], are often applicable to small and moderate-sized problems, but are often too expensive to be practical for large sparse problems. Moreover, direct methods may be susceptible to numerical instability and the loss of solution accuracy, as evidenced by the findings of [5].

In contrast to most iterative methods, direct methods are typically easier to implement. This is particularly true when a coefficient matrix has a specific structure, such as that shown in (2). In such cases, it is possible to exploit the structure of the matrix to create a fast solver that is as easy to implement as iterative methods but faster than classical direct methods, such as Gauss elimination, LU, and so on.

This paper is organized as follows: The following section presents some fast algorithms for solving (1) with the coefficient matrix T being a tridiagonal Toeplitz matrix, analogous to the matrix in (2). Additionally, B is an $n \times m$ matrix containing multiple right-hand side vectors. The efficacy of our algorithm is demonstrated through numerical tests in Section 3, after which a brief conclusion is drawn in Section 4.

## 2 Fast Algorithm

If we consider (1) as follow:

$$TXI_m = B \tag{3}$$

with use of Kronecker product and **vec** operator, we have

$$\underbrace{(I_m \otimes T)}_{A} x = b \tag{4}$$

where

$$x = vec(X), \quad b = vec(B)$$

In above, the matrix $A$ hase the form:

$$A = \begin{bmatrix} T & & & \\ & T & & \\ & & \ddots & \\ & & & T \end{bmatrix} \tag{5}$$

but $A$ itself is not full teridiagonal Toeplitz.

Let

$$E_1 = \begin{bmatrix} 0 & \cdots & 0 \\ \vdots & \ddots & \vdots \\ t_3 & \cdots & 0 \end{bmatrix}, E_2 = \begin{bmatrix} 0 & \cdots & t_2 \\ \vdots & \ddots & \vdots \\ 0 & \cdots & 0 \end{bmatrix} \tag{6}$$

and

$$J = \begin{bmatrix} 0 & 1 & & & \\ 0 & & 1 & & \\ \vdots & & & \ddots & \\ 0 & & & & 1 \\ 0 & 0 & \cdots & 0 & 0 \end{bmatrix}$$

,

Matrix

$$\hat{A} = A + (J \otimes E_1) + (J^T \otimes E_2) \tag{7}$$

is a full tridiagonal Toeplitz matrix.

We can use another representation for the matrix $\hat{A}$ as follow:

$$\hat{A} = A + \sum_{j=1}^{m-1} \alpha_{jm} e_{jm+1}^T + \beta_{jm+1} e_{jm}^T \tag{8}$$

where $e_j \in \mathbb{R}^{mn \times mn}$ denots the $jth$ colmn of Identity matrix $I_{mn}$, and

$$\alpha_j = (0, \cdots, 0, \underbrace{t_3}_{jth}, 0, \cdots, 0)^T \in \mathbb{R}^{mn \times mn}$$

$$\beta_j = (0, \cdots, 0, \underbrace{t_2}_{jth}, 0, \cdots, 0)^T \in \mathbb{R}^{mn \times mn}$$

Now, with (8), we can rewrite the (4) as follow:

$$(\hat{A} - \sum_{j=1}^{m-1} \alpha_{jm} e_{jm+1}^T + \beta_{jm+1} e_{jm}^T)x = b \tag{9}$$

$\hat{A}$ is a full tridiagonal Toeplitz matrix and solving systems of equations with $\hat{A}$ as their coefficient matrix can be done fast with help of a fast algorithm. In this paper we use algorithm introduced in [7]. This fast algorithm uses $2 \times 2$ Block LU-factorization and backward sunstitution to solve tridiagonal Toeplitz systems, and since our original coefficient matrix, $A$, is not full tridiagonal Toeplitz, algorithm will not work correctly. By looking closely to elements of matrix $A$, it's easy to see that there are $2m - 2$ zeros on the upper and lower diagonal of matrix $A$ and the Idea of this paper is to add an approprate matrix to the $A$ and make it full tridiagonal Toeplitz so use of the fast algorithm in [7] could result more reliable solution.

We have

$$(I - \sum_{j=1}^{m-1} \alpha_{jm} e_{jm+1}^T \hat{A}^{-1} + \beta_{jm+1} e_{jm}^T \hat{A}^{-1})\hat{A}x = b$$

Let

$$\hat{A}x = \phi, \tag{10}$$

Then

$$(I - \sum_{j=1}^{m-1} \alpha_{jm} e_{jm+1}^T \hat{A}^{-1} + \beta_{jm+1} e_{jm}^T \hat{A}^{-1})\phi = b \tag{11}$$

Now, for any $i$, let $y_i$ be the solution of the $\hat{A}^T y_i = e_i$. Then from (11) we have

$$(I - \sum_{j=1}^{m-1} \alpha_{jm} y_{jn+1}^T + \beta_{jm+1} y_{jn}^T)\phi = b \tag{12}$$

**Lemma :** ([8], p.563) If A is an $n \times n$ matrix and $U_k$ and $V_k$, $(k = 1, ..., N)$ are $n \times m$ matrices then

$$\left(A + \sum_{k=1}^{N} U_k V_k^T\right)^{-1} = A^{-1} - A^{-1}[U_1, U_2, ..., U_N]M^{-1}[V_1^T, V_2^T, ..., V_N^T]A^{-1}$$

Where $M$ is $Nm \times Nm$ matrix given by

$$M = \begin{bmatrix} I_{m\times m} + V_1^T A^{-1} U_1 & V_1^T A^{-1} U_2 & \cdots & V_1^T A^{-1} U_N \\ V_2^T A^{-1} U_1 & I_{m\times m} + V_2^T A^{-1} U_2 & \cdots & V_2^T A^{-1} U_N \\ \vdots & & & \vdots \\ V_N^T A^{-1} U_1 & V_N^T A^{-1} U_2 & & I_{m\times m} + V_N^T A^{-1} U_N \end{bmatrix}.$$

**Corollary :**([8], p.563) Let $A = I$ the $n \times n$ identity matrix. Then

$$\left(I + \sum_{k=1}^{N} U_k V_k^T\right)^{-1} = I - [U_1, U_2, ..., U_N]M^{-1}[V_1^T, V_2^T, ..., V_N^T]$$

Where $M$ is $Nm \times Nm$ matrix given by

$$M = \begin{bmatrix} I_{m\times m} + V_1^T U_1 & V_1^T U_2 & \cdots & V_1^T U_N \\ V_2^T U_1 & I_{m\times m} + V_2^T U_2 & \cdots & V_2^T U_N \\ \vdots & & & \vdots \\ V_N^T U_1 & V_N^T U_2 & & I_{m\times m} + V_N^T U_N \end{bmatrix}.$$

If we set

$$Q = (I - \sum_{j=1}^{m-1} \alpha_{jm} y_{jn+1}^T + \beta_{jm+1} y_{jn}^T)$$

with use of above Corollary we can compute $Q^{-1}$ and from (12) we have

$$\phi = Q^{-1}b \tag{13}$$

Then by replacing $\phi$ in (10)

$$\hat{A}x = Q^{-1}b \tag{14}$$

$\hat{A}$ in (14) is tridiagonal Toeplitz matrix and fast solver in [7] can be utilized to solve (14). Summarizing above analysis gives us the following algorithm:

**Algorithm 1** Fast Solver For Tridiagonal Toeplitz Systems with Multiple RHS

**Require:** Tridiagonal Toeplitz Matrix $T$, RHS Matrix $B$,

1: Determine $\alpha_j$ and $\beta_j$ as (8).

2: Use the Fast Solver in [7] to solve $2m-2$ systems $\hat{A}^T y_i = e_i$ correspond to the appropriate indexs of $e_i$ in (11).

3: Put

$$Q = (I - \sum_{j=1}^{m-1} \alpha_{jm} y_{jn+1}^T + \beta_{jm+1} y_{jn}^T)$$

and compute $Q^{-1}$ with the help of corollary of Sherman-Morrison-Woodbury Formula.

4: Use the Fast Solver in [7] to solve tridiagonal Toeplitz system (14)

5: **return** x.

## 3 Numerical Examples

In this section we use some examples to show effectiveness of proposed Algorithm. All the numerical tests were done on an ASUS laptop PC with AMD A12 CPU, 8Gb RAM and by Matlab R2016(b) with a machine precision of 1016. For convenience, throughout our numerical experiments, we denote by Relative error = $||B - TX||/||B||$, the relative residual error, and computing time (in seconds). In all tables, the Time is the average value of computing times required by performing the corresponding algorithm 10 times, the right-hand side matrix B is taken to be $B = ones(n, m)$.

**Example 1:** For our first example we consider the matrix T be a tridiagonal "Grcar" matrix which is in following form:

$$T = \begin{bmatrix} 1 & 1 & & & \\ -1 & 1 & 1 & & \\ & \ddots & \ddots & \ddots & \\ & & -1 & 1 & 1 \\ & & & -1 & 1 \end{bmatrix}$$

Table(1), shows the performance of Algorithm(1) for different matrix size and different number of right-hand sides for example 1.

Table 1: Reslts for example 1 :

| n | m | Relative Error | Time | n | m | Relative Error | Time |
|---|---|---|---|---|---|---|---|
| 10 | 2 | 6.4370 e-13 | 3.6972 e-4 | 20 | 2 | 2.7244 e-9 | 9.8775 e-4 |
| | 3 | 3.3816 e-11 | 9.0149 e-4 | | 3 | 2.0428 e-5 | 0.0019 |
| | 4 | 4.3365 e-9 | 0.0020 | | 4 | 0.8686 | 0.0031 |
| | 5 | 3.3256 e-7 | 0.0036 | | 5 | 5.9842 e+3 | 0.0072 |
| | 6 | 6.3730 e-5 | 0.0047 | 30 | 2 | 3.6752 e-5 | 0.0012 |
| | 7 | 0.0058 | 0.0057 | | 3 | 894428 | 0.0025 |
| | 8 | 0.9774 | 0.0059 | | 4 | 5.6572 e+7 | 0.0069 |
| | 9 | 166.6069 | 0.0099 | 40 | 2 | 1.2796 | 0.0015 |
| | 10 | 3.3036 e+3 | 0.0160 | | 3 | 4.1251 e+7 | 0.0065 |

**Example 2:** In this example we consider T to be a real nonsymmetric Tridiagonal Toeplitz matrix created by $f(\theta) = e^{i\theta} + 2e^{-i\theta}$ as follow:

$$T = \begin{bmatrix} 0 & 2 & & \\ 1 & 0 & \ddots & \\ & \ddots & \ddots & 2 \\ & & 1 & 0 \end{bmatrix}$$

Table(2) shows the performance of Algorithm(1) for different matrix size and different number of right-hand sides for example 2.

Table 2: Results for Example 2:

| n | m | Relative Error | Time |
|---|---|---|---|
| 10 | 2 | 0 | 5.6249 e-4 |
| | 4 | 0 | 0.0032 |
| | 8 | 0 | 0.0090 |
| | 10 | 0 | 0.0139 |
| 30 | 2 | 0 | 0.0012 |
| | 4 | 0.2236 | 0.0093 |
| | 8 | 0.5105 | 0.0715 |
| | 10 | 0.5553 | 0.1531 |
| 50 | 2 | 0 | 0.0027 |
| | 4 | 0.4695 | 0.0286 |
| | 8 | 0.5980 | 0.2836 |
| | 10 | 0.6212 | 0.5610 |

## 4 Conclusion

In this work, we introduce a fast algorithm for solving Tridiagonal Toeplitz linear systems with Multiple right-hand Sides. First we added $2m - 2$ rank one matrix to matrix $A$ and filled $2m - 2$ zeros on its upper and lower diagonal, and turned it into a full tridiagonal Toeplitz matrix. Then with help of Sherman-Morrison-Woodbury formula, we rearranged our problem into a problem consist of full tridiagonal Toeplitz matrix $\hat{A}$. Then we used thae Tridiagonal Toeplitz Fast Solver algorithm intriduced in [7] to solve our problem fast.

**Future Works :** Numerical results in section 3 show that algorithm(1), like most direct methods, is not suitable for large-scale problems. Therefore, future work(s) would be on finding a reliable block method to extend the fast algorithm to large-scale problems.

## References

[1] R. Chan and X.-J. Jin, An introduction to iterative Toeplitz solvers. *SIAM, Philadelphia, PA, USA, 2007.*

[2] G. Golub and C. Van Loan, Matrix Computations, *Johns Hopkins University Press, Baltimore and London, 3rd edition, 1996.*

[3] Y. Saad, Iterative methods for sparse Linear systems, *2nd edition, SIAM, Philadelphia, PA, USA, 2003.*

[4] Z.-Y. Liu, X.-R. Qin, N.-C. Wu and Y.-L. Zhang, The shifted classical circulant and skew circulant splitting iterative methods for Toeplitz matrices, *Canad. Math. Bull., 60 (2017): 807-815.*

[5] R. Chan and M. K. Ng, Conjugate gradient methods for Toeplitz systems, *SIAM Rev. 38 (1996): 427-482.*

[6] S.S. Nemani, L.E. Garey, An efficient method for solving second order boundary value problems with two point boundary conditions, *Int. J. Comput. Math. 79 (9) (2002) 10011008.*

[7] ZHONGYUN LIU , SHAN LI , YI YIN , AND YULIN ZHANG FAST SOLVERS FOR TRIDIAGONAL TOEPLITZ LINEAR SYSTEMS, *09 November 2020*

[8] Milan Batista, Abdel Rahman A, Ibrahim Karawia The use of the ShermanMorrisonWoodbury formula to solve cyclic block tri-diagonal and cyclic block penta-diagonal linear systems of equations *Appl. Math. Comput. 210 (2009) 558-563*